\documentclass[11pt]{article}

\usepackage[psamsfonts]{amssymb}
\usepackage{amsfonts,euscript,wrapfig,url}
\usepackage{amsmath,amsthm,graphicx,enumitem}
\usepackage{booktabs,color,pgf,tikz,pgfbaseimage}
\usepackage{mathrsfs}
\usetikzlibrary{arrows}

\title{Self-Similar Dendrites Generated by Polygonal Systems in the Plane}

\author{Mary Samuel \and Andrey Tetenov \footnote{Supported by Russian Foundation of Basic Research project 16-01-00414}\and Dmitry Vaulin }

\begin{document}
\newcommand{\rr}{\mathbb{R}}
\newcommand {\nn} {\mathbb{N}}
\newcommand {\zz} {\mathbb{Z}}
\newcommand {\bbc} {\mathbb{C}}
\newcommand {\rd} {\mathbb{R}^d}
\newcommand {\rpo}{\mathbb{R}_+^1}

 \newcommand {\al} {\alpha}
\newcommand {\be} {\beta}
\newcommand {\da} {\delta}
\newcommand {\Da} {\Delta}
\newcommand {\ga} {\gamma}
\newcommand {\Ga} {\Gamma}
\newcommand {\la} {\lambda}
\newcommand {\La} {\Lambda}
\newcommand{\om}{\omega}
\newcommand{\Om}{\Omega}
\newcommand {\sa} {\sigma}
\newcommand {\Sa} {\Sigma}
\newcommand {\te} {\theta}
\newcommand {\fy} {\varphi}
\newcommand {\Fy} {\varPhi}

\newcommand {\rh} {\varrho}

\newcommand{\ep}{\varepsilon}

\newcommand{\VEC}{\overrightarrow}
\newcommand {\ra} {\rightarrow}
\newcommand{\IN}{{\subset}}
\newcommand{\ov}{{\overline}}
\newcommand{\NI}{{\supset}}
\newcommand \dd  {\partial}
\newcommand {\mmm}{{\setminus}}
\newcommand{\probel}{\vspace{.5cm}}
\newcommand{\8}{{\infty}}
\newcommand{\io}{{I^\infty}}
\newcommand{\0}{{\varnothing}}
\newcommand{\vse}{$\blacksquare$}

\newcommand {\bfep} {{{\bar \varepsilon}}}
\newcommand {\Dl} {\Delta}
\newcommand{\vA}{{\vec {A}}}
\newcommand{\vB}{{\vec {B}}}
\newcommand{\vF}{{\vec {F}}}
\newcommand{\vf}{{\vec {f}}}
\newcommand{\vh}{{\vec {h}}}
\newcommand{\vJ}{{\vec {J}}}
\newcommand{\vK}{{\vec {K}}}
\newcommand{\vP}{{\vec {P}}}
\newcommand{\vX}{{\vec {X}}}
\newcommand{\vY}{{\vec {Y}}}
\newcommand{\vZ}{{\vec {Z}}}
\newcommand{\vx}{{\vec {x}}}
\newcommand{\va}{{\vec {a}}}
\newcommand{\vga}{{\vec {\gamma}}}

\newcommand{\hf}{{\hat {f}}}
\newcommand{\hg}{{\hat {g}}}

\newcommand{\bj}{{\bf {j}}}

\newcommand{\bi}{{\bf {i}}}
\newcommand{\bk}{{\bf {k}}}
\newcommand{\bu}{{\bf {u}}}

\newcommand{\bX}{{\bf {X}}}
\newcommand{\Ep}{{\mathfrak p}}
\newcommand{\Eq}{{\mathfrak q}}
\newcommand{\eB}{{\EuScript B}}
\newcommand{\wP}{{\widetilde P}}
\newcommand{\eU}{{\EuScript U}}
\newcommand{\eS}{{\EuScript S}}
\newcommand{\eH}{{\EuScript H}}
\newcommand{\eC}{{\EuScript C}}
\newcommand{\eP}{{\EuScript P}}
\newcommand{\eT}{{\EuScript T}}
\newcommand{\eG}{{\EuScript G}}
\newcommand{\eK}{{\EuScript K}}
\newcommand{\eF}{{\EuScript F}}
\newcommand{\eZ}{{\EuScript Z}}
\newcommand{\eL}{{\EuScript L}}
\newcommand{\eD}{{\EuScript D}}
\newcommand{\E}{{\EuScript E}}
\def \diam {\mathop{\rm diam}\nolimits}
\def \sup {\mathop{\rm sup}\nolimits}
\def \fix {\mathop{\rm fix}\nolimits}
\def \Lip {\mathop{\rm Lip}\nolimits}
\def \min {\mathop{\rm min}\nolimits}

\newcommand{\red}{\textcolor{red}}

\newtheorem{thm}{\bf Theorem}
 \newtheorem{cor}[thm]{\bf Corollary}
 \newtheorem{lem}[thm]{\bf Lemma}
 \newtheorem{prop}[thm]{\bf Proposition}
 \newtheorem{dfn}[thm]{\bf Definition}
 \newcommand{\rmk}{{\bf {Remark.}}}

\newcommand{\dok}{{\bf{Proof:  }}}

\maketitle

\bigskip

{\bf Abstract}: We consider a method of construction of self-similar dendrites on a plane and establish  main topological and metric properties of resulting class of dendrites.

\smallskip

\section{Introduction}

The study of dendrites occupies a significant place in general topology \cite{Kur,Nad,Whyb}. One can refer to the paper \cite{Char} of J.Charatonik and W.Charatonik for   exhaustive overview covering more than 75-year research in this area. At the same time, in  the theory of self-similar sets there are only individual  attempts to work out some approaches to self-similar dendrites in certain situations. In 1985, Hata \cite{Hata} studied  topological properties of  attractor $K$ of  a system $\eS$ of weak 
contractions in a complete metric space and showed that if  $K$ is a dendrite then it has infinite set of end points. Jun Kigami in his work \cite{Kig95} applied the methods of harmonic calculus on fractals to dendrites; on a way to this he developed effective approaches to the study of  structure of self-similar dendrites. D.Croydon in his thesis \cite{C} obtained heat kernel estimates for continuum random tree and for certain family of p.c.f. random dendrites on the plane. D.Dumitru and A.Mihail \cite{DM} made an attempt to get a sufficient condition for a self-similar set to be a dendrite in terms of sequences of  intersection graphs for the refinements of the system $\eS$. We need also mention very useful ideas in \cite{BR} and examples in  \cite[Fig.VII.200]{Bar}.

There are several questions arising in the study of self-similar dendrites.
What kind of topological restrictions characterise the class of dendrites generated by systems of similarities in $\rr^d$? What are the explicit construction algorithms for self-similar dendrites? What are the metric and analytic properties of morphisms of self-similar structures on dendrites?

The aim of our work is to make clear basic topological and metric properties of self-similar dendrites in the most simple  settings. For that reason we consider systems of similarities in the plane, which we call polygonal tree systems (Definition \ref{pts}). We show that the attractor $K$ of such system $\eS$ is a  dendrite (Theorem \ref{main}), that, by the  construction, each such system  $\eS$  satisfies open set condition, one-point intersection property and   is post-critically finite (Proposition \ref{obvious}, Corollary \ref{pcf}); for the dendrite $K$ we define its main tree (Definition \ref{defmt}) and show that each cut point of $K$ lies in some image $S_\bj(\hat\ga)$ of the main tree (Theorem \ref{order}) and get the upper bound for the order of   ramification points of $K$, depending only on the initial polygon $P$ of the system $\eS$. We show that the dendrite $K$ is a continuum with bounded turning in the sense of P.Tukia  (Theorem \ref{c-BT}). Finally,  we show that each combinatorial equivalence of polygonal tree systems $\eS$,$\eS'$ defines  unique  homeomorphism $\fy: K\to K'$, compatible  with $\eS$ and $\eS'$ and prove H\"older continuity of $\fy$ and $\fy^{-1}$ (Theorem \ref{morphism}).

\subsection{Preliminaries} 

{\bf Dendrites.}  A {\em dendrite} is a locally connected continuum containing no simple closed curve. 

We shall use the notion of {\it order of a point} in the sense of Menger-Urysohn (see \cite[Vol.2, \S 51, p.274]{Kur}) and we denote by $Ord (p, X)$ the order of the 
continuum $X$ at a point $p \in X$. 
{Points of order 1 in a continuum $X$ are called {\em end points} of $X$;   the set of all end points of $X$   will
be  denoted by $EP(X)$. A point $p$ of a continuum $X$ is called a {\em cut point} of $X$ provided that $X \setminus \{p\}$ is
not connected; the set of all cut points of $X$ will be denoted by $CP(X)$. 
Points of order at least 3 are called {\em ramification points} of $X$; the  set of all ramification points of $X$ will be denoted by $RP(X)$. }

We will  use the following statements selected from  \cite[Theorem 1.1]{Char}:
\begin{thm} For a continuum $X$ the following conditions are equivalent:
\begin{enumerate}[label=(\alph*),noitemsep]
\item $X$ is dendrite;
\item every two distinct points of $X$ are separated by a third point;
\item each point of $X$ is either a cut point or an end point of $X$;
\item each nondegenerate subcontinuum of $X$ contains uncountably many cut points of $X$.
\item for each point $p \in X$ the number of components of the set $X \setminus \{p\} = ord (p, X)$ whenever either of these is finite;
\item the intersection of every two connected subsets of X is connected;
\item X is locally connected and uniquely arcwise connected.
\end{enumerate}
\end{thm}

\bigskip
{\bf Self-similar sets.}
Let $(X, d)$ be a complete metric space. 
A mapping\\ $F: X \to X$ is a contraction if $\Lip F < 1$. 
The mapping $S: X \ra X$ is called a similarity if \begin{equation} d(S(x), S(y)) = r d(x, y) \end{equation} for all $x, y\in X$ and some fixed r. 

\begin{dfn} 
Let $\eS=\{S_1, S_2, \ldots, S_m\}$ be a system of (injective) contraction maps on the complete metric space $(X, d)$.
 A nonempty compact set $K\IN X$ is said to be invariant with respect to $\eS$, if $K = \bigcup \limits_{i = 1}^m S_i (K)$. \end{dfn}
 
 We also call the   subset $K \IN X$ self-similar with respect to $\eS$. Throughout the whole paper, the maps $S_i\in \eS$ are supposed to be  similarities and the set $X$ to be $\rr^2$.\\
 We denote $I=\{1,2,\ldots,m\}$,  $I^*=\bigcup_{n=1}^\8I^n$ is the set of all finite $I$-tuples $\bj=j_1j_2...j_n$, 
  $I^{\8}=\{{\bf \al}=\al_1\al_2\ldots,\ \ \al_i\in I\}$ is the index  space and 
 $\pi:I^{\8}\rightarrow K$ is the address map.\\ 
 As usual for any $\bj\in I^*$, we  write $S_\bj=S_{j_1j_2...j_n}=S_{j_1}S_{j_2}...S_{j_n}$ and for some set $A\IN X$ we  often denote $S_\bj(A)$ by $A_\bj$.

 \begin{dfn}The system ${\eS}$ satisfies the {\em open set condition} (OSC) if there exists a non-empty open set $O \IN X$ such that $S_i (O), \{1 \le i\le m\}$ are pairwise disjoint and all contained in $O$.\end{dfn}
 
We say the self-similar set $K$ defined  by  the  system $\eS$  satisfies the one-point intersection property if for  any $i\neq j$,  $S_i(K) \bigcap S_j(K)$ is not more than one point. 

The union $\eC$ of all  $S_i(K)\cap S_j(K)$, $i,j \in I, i\neq j$ is called  the critical set of  the system $\eS$. 
The post-critical set $\eP$ of the system $\eS$ is the set of all
 $\alpha\in I^{\8}$ 
such that for some ${\bf j}\in I^*$,
 $S_ {\bf j}(\pi(\al))\in\eC$. \cite{Kig}
\smallskip

{\bf Kigami's theorem.} We use  the following  convenient criterion of  connectedness of  the  attractor  of a system $\eS$    \cite{Kig}:

\begin{dfn} \label{conn}
Let $ \{S_i(K)\}_{i \in I}, \{I = 1, 2, \ldots, m\}$  be a family of non-empty subsets of $X$. The family $ \{S_i(K)\}_{i \in I}$ is said to be connected if for every $i, j \in I$ there exists $\{i_0, i_1, \ldots, i_n\} \IN I$ such that $i_0 = i, i_n = j$ and \\$S_{i_k}(K)\bigcap S_{i_{k + 1}}(K) \ne \0$ for every $k = 0, 1, \ldots, {n-1}$.\end{dfn}

\begin{thm}\label{Kigami} Let $(X, d)$ be a complete metric space on which a finite number of contractions $S_i: X \rightarrow X$ are defined such that the self-similar set corresponding to the system of maps be  $K = \bigcup \limits_{i =1}^m S_i (K)$. Then the following statements are equivalent:\\ (1) The family $\{S_i (K)\}_{i = 1}^m$ is connected.\\(2) $K$ is arcwise connected.
\\
(3) $K$ is connected. \end{thm}

\bigskip

{\bf Zippers and multizippers.} The simplest  way to construct  a self-similar  curve is  to  take a polygonal line and then replace each of its segments by a smaller copy of the same polygonal line; this construction is called   zipper and was studied  by Aseev, Tetenov and  Kravchenko \cite{ATK}.

\begin{dfn} Let $X$  be a complete metric space. A system $\eS = \{S_1, \ldots, S_m\}$
of contraction mappings of $X$ to itself is called a {\em zipper} with vertices $\{z_0, \ldots, z_m\}$
and signature $\ep = (\ep_1, \ldots, \ep_m), \ep_i \in\{0,1\}$, if for    $ i = 1\ldots m, S_i (z_0) = z_{i-1+\ep_i}$ and $S_i (z_m) = z_{i-{\ep_i}}$.\end{dfn}

More  general approach  for  building self-similar  curves and  continua is provided  by  a graph-directed version of zipper construction \cite{Tet06}:

\begin{dfn} Let $\{X_u, u \in V\}$ be a system of spaces, all isomorphic to ${\rr}^d$.
 For each $X_u$ let a finite array of points be given $\{x_0^{(u)}, \ldots, x_{m_u}^{(u)}\}$.
 Suppose for each $u\in V$ and $0\le k \le m_u$ we have some $v (u, k) \in V$  and  $\ep(u,k)\in\{0,1\} $   and a map ${\eS}_k^{(u)}: X_v \ra X_u$ such that\\
 $S^{(u)}_k (x_0^{(v)}) = x_{k-1}^{(u)}$  or $x_k^{(u)}$  and  
$S_k^{(u)} (x_{m_v}^{(v)}) = x_{k}^{(u)}$  or  $  x_{k-1}^{(u)}$,
depending on the signature $\ep (u, r)$.\\ The graph directed iterated function system (IFS) defined by the maps $S_k^{(u)}$ is called a {\em multizipper} ${\eZ}.$\end{dfn}

The attractor of multizipper ${\eZ}$ is a system of  connected and  arcwise connected  compact  sets $K_u\IN X_u$ satisfying  the system of equations
$$ K_u=\bigcup\limits_{k=1}^{m_u}S_k^{(u)}(K_{v(u,k)}),\qquad u\in V$$
 We call the  sets $K_u$  the components of  the attractor of ${\eZ}$.
   
 %These  components are the H\"older images of $[0,1]$:
%\begin{prop} Let $\{X_u, u \in V\}$ be a system of spaces and $\eZ =\{S_e\}$ be a multizipper with  vertices $\{x_i^{(u)}, u \in V, i = 1, 2, \ldots m_u\}$ and signature $\overline \ep = \{(v (u, k), \ep (u, k))\}$. For a linear multizipper $\eL$ with the same signature and any choice of node points $O = t_0^{(u)} < t_1^{(u)} < \ldots t_{m_u}^{(u)} = 1$ on intervals $I_u = [0, 1]$ there is unique array of surjective continuous maps $$f_u : I_u \ra K_u$$ such that for any $u \in V$ and for any $k \in \{1, \ldots m_u\}$.\\ $$f_u (t_k ^{(u)}) = x_k^{ (u)}$$ and $$S_k^{(u)} \circ f_{v(u, k)} = f_u \circ T_{v (u,k)}$$ All $f_u$ being Holder continuous. \end{prop}

The components $K_u$  of the  attractor of $\eZ$ are  Jordan  arcs if  the following  conditions  are  satisfied:
 
\begin{thm} \label{Jor}
 Let ${\eZ}_0 = \{S_k^{(u)}\}$ be a multizipper with node points $x_k ^{(u)}$ and a signature $\ep = \{(v (u, k), \ep (u, k)), u \in V, k = 1, \ldots, m_u\}$. If for any $u \in V$
and any  $ i, j \in \{1, 2, \ldots, m_u\}$, the set $K_{(u, i)} \cap K_{(u, j)} = \0$ if $|i- j| > 1$ and is a singleton if $|i -j| = 1$, then any linear parametrization $\{f_u: I_u \ra K_u\}$ is a homeomorphism and each $K_u$ is a Jordan arc with endpoints $x^{(u)}_0,\   x_m ^{(u)}.$\end{thm}
\section{Polygonal tree  systems.}
Let $P$ be a convex polygon in ${\rr}^2$ and $A_1, \ldots, A_n$ be its vertices.\\
Consider a system of contracting similarities $\eS = \{S_1, \ldots, S_m\}$, which possesses  the following properties:\\
{\bf(D1)}\ \  For any $k = 1, \ldots, m$, the set $P_k = S_k (P)$ is contained in $P$;
\\
{\bf(D2)}\ \   For any $i\neq j,\ \   i, j = 1, \ldots, m, P_i \bigcap P_j$ is either empty or is a common vertex of $P_i$ and $P_j$;\\
  {\bf(D3)}\ \  For any vertex $A_k$ there is  the map $S_i\in\eS$ such that $P_i\ni A_k$;\\
{\bf(D4)}\ \   {\em The set    ${\wP} = \bigcup \limits_{i = 1}^m P_i$ is contractible.}

\begin{dfn}\label{pts}
The system $(P,\eS)$ satisfying the conditions {\bf D1-D4} is called  a polygonal tree system associated with the polygon $P$.
\end{dfn}

Some properties of the attractor $K$ of a polygonal tree system $\eS$ follow  directly  from its  definition:
\begin{prop}\label{obvious} Let $\eS$ be a polygonal tree system associated  with a polygon $P$ and let $K$ be its attractor. Then \quad
(i) $\eS$ satisfies open set  condition; \quad (ii) $\eS$ satisfies one point intersection property. 
\end{prop}

\dok (i) Since for any $i,j=1,\ldots,m$, $P_i\IN P$ and $\dot{P_i}\cap\dot{P_j}=\0$ for $i\neq j$, $\dot P$ can be taken for the open set; (ii) follows directly from {\bf(D2)}\vse \\

Thus, to define a polygonal tree  system we specify  a polygon $P$, a system of its subpolygons $P_i$ and the similarities $S_i$, sending $P$ to $P_i$. Along with each polygonal tree system ${\eS}$ we consider the set ${\wP} = \bigcup \limits_{i = 1}^{m} P_i$ and  the Hutchinson operator $H_S(A)=\bigcup\limits_{i=1}^m S_i(A)$ which sends $P$ to $\wP$.

\newpage
\noindent{\bf Example 2.1. Hata's tree-like set.}\bigskip\\
\begin{tikzpicture}\node[anchor=south west,inner sep=0] at (0,0) {
   \includegraphics[width=.52\textwidth]{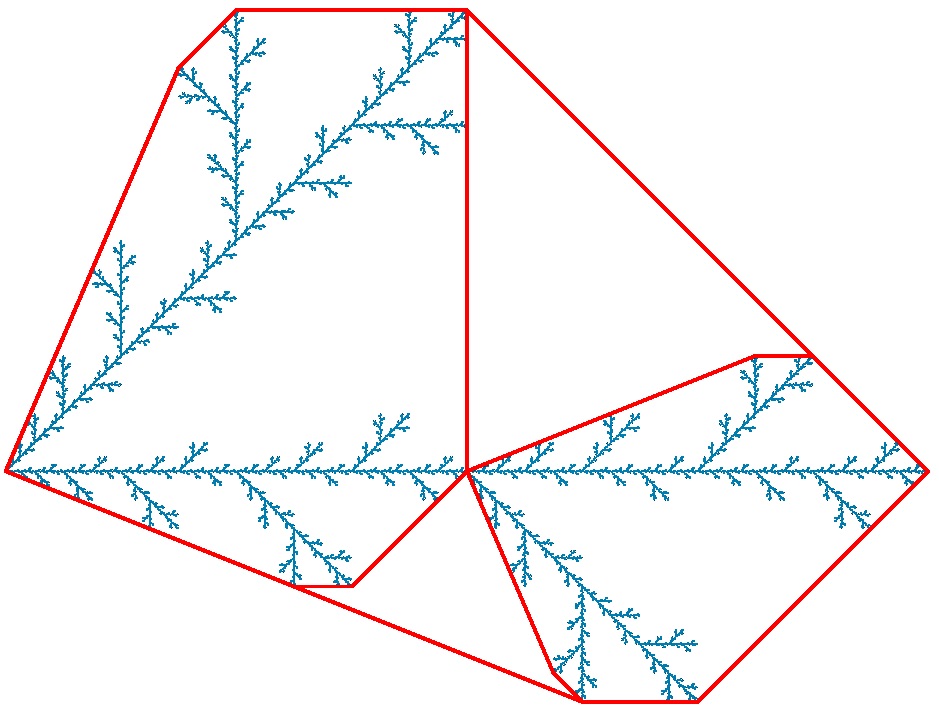}};
\draw (-1.9,1.8) node[anchor=north west] {$\Large A_1=S_1(A_1)$};
\draw (2.6,0.12) node[anchor=north west] {$\Large A_2=S_2(A_6)$};
\draw (4.9,.21) node[anchor=north west] {$\Large A_3=S_2(A_5)$};
\draw (6.4,2) node[anchor=north west] {$ A_4=S_2(A_4)$};
\draw (3.47,5.05) node[anchor=north west] { $A_5=S_1(A_4)$};
\draw (1.2,5.5) node[anchor=north west] {$\Large A_6=S_1(A_3)$};
\draw (-1,4.8) node[anchor=north west] {$\Large A_7=S_1(A_2)$};
\draw (1.92,2.8) node[anchor=north west] {\Large $P_1$};
\draw (4.6,1.2) node[anchor=north west] {\Large $P_2$};
\end{tikzpicture}
\medskip\\
{\small
 Hata's  tree-like set \cite{C, Hata, Kig} is  the attractor of a polygonal tree system. The polygon $P$ for the  set has 7 vertices. The maps are $S_1(z)=(1+i)\bar z/2$,\quad $S_2(z)= (\bar z+1)/2$.}\bigskip \\ 
 {\bf Example 2.2.}
\begin{center}
\begin{tikzpicture}[scale=5.5]
\node[anchor=south west,inner sep=0] at (-.02,-.02) {
   \includegraphics[width=.45\textwidth]{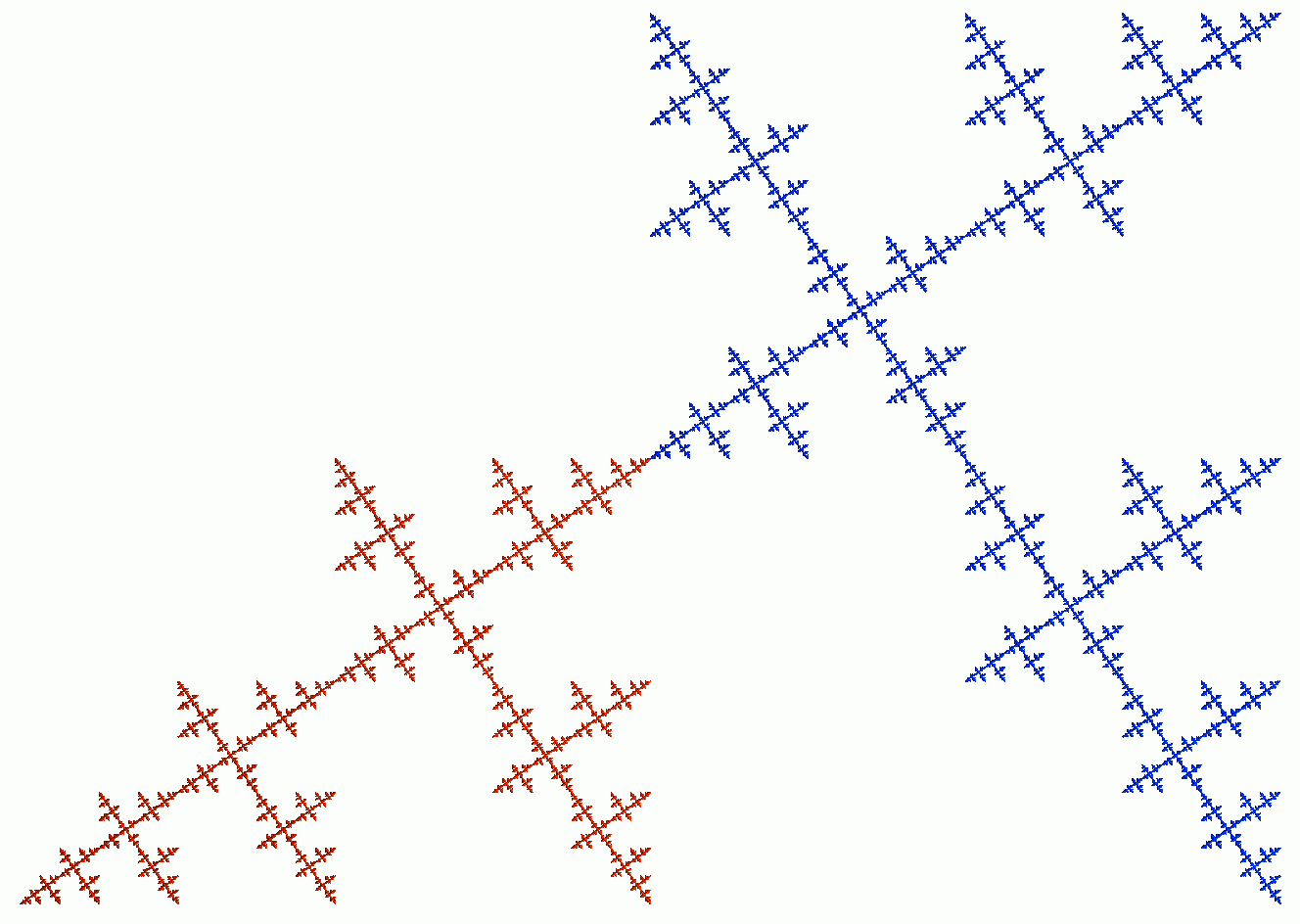}};

\fill[color=green,fill=green,fill opacity=0.25] (0.,0.) -- (0.25,0.3536) -- (0.5,0.3536) -- (0.5,0.) -- cycle;
\fill[color=green,fill=blue,fill opacity=0.25] (0.5,0.707) -- (0.5,0.3536) -- (1.,0.) -- (1.,0.707) -- cycle;
\draw [line width=1.2pt,color=red] (0.,0.)-- (0.5,0.707)-- (1.,0.707)-- (1.,0.)-- (0.,0.);
\draw (0.24,0.2) node[anchor=north west] {\Large $P_1$};
\draw (0.72,0.5) node[anchor=north west] {\Large $P_2$};
\draw (-0,0.06) node[anchor=north east] {$\Large A_1=S_1(A_1)$};
\draw (1.,0.06) node[anchor=north west] {$\Large A_2=S_2(A_1)$};
 
\draw (0.5,0.77) node[anchor=north east] {$\Large A_4=S_2(A_3)$};
\draw (1,0.77) node[anchor=north west] {$\Large A_3=S_2(A_2)$};
\end{tikzpicture}\\
{\small  A polygonal tree system   $(P,\eS)$,   $S_1 (z) = z/2$ and  $S_2 (z) = {i z}/\sqrt{2}+1$ defines a dendrite from R.Zeller's thesis   \cite[Ch.1, p.18]{Zell}.}
\end{center}

\bigskip
Composition  of  two  Hutchinson operators corresponding  to  two polygonal tree  systems  associated  with the same polygon $P$  is  also  an operator of  the  same  type:

\begin{lem}\label{trans}   Let $(P,\eS)$ and ${(P,\eS')}$ be polygonal tree systems of similarities associated with $P$. Then the system ${\eS}'' = \{S_i \circ S_j', S_i \in {\eS}, S_j \in {\eS}'\}$ is a polygonal tree system of similarities associated with  $P$.\end{lem}
\noindent{\dok} {\bf(D1)}\ \ is obvious because $S_i \circ S_j'(P) \IN S_i (P) \IN P$.\\
{\bf(D2)}\ \ Let $Q_1 =S_{i_1}\circ S_{j_1}' (P)$ and $Q_2 = S_{i_2}\circ S_{j_2}' (P)$ be two polygons in ${\eS}''$ and consider their intersection:\\
if $i_1 \ne i_2$,  $Q_1 \bigcap Q_2 \IN P_{i_1} \bigcap P_{i_2}$, where the left-hand side intersection contains  at most one point.\\
if $i_1 = i_2, Q_1 \bigcap Q_2 = S_{i_1} (P_{j_1}' \bigcap P_{j_2}')$ which is either empty or a one-point set, containing $S_{i_1} (A')$ where $A'$ is   a common vertex of $P_{j_1}'$ and $P_{j_2}'$.\\
{\bf(D3)}\ \   holds because for any vertex $A_k$, the similarity $S_k \circ S_k'$ is the unique similarity in ${\eS}''$, fixing the point $A_k$.
\\
{\bf(D4)}\ \ The sets $\wP=\bigcup \limits_{i = 1}^m P_i$ and $\wP'=\bigcup
\limits_{i = 1}^{m'} P_i'$ are strong deformation retracts
of the polygon $P$, both containing the vertices $A_1, \ldots A_n$ of
$P$.
Let ${\fy}' (X, t): P \times [0, 1] \ra P$ be a deformation
retraction from $P$ to $\mathop{\bigcup}\limits_{i = 1}^{m'} P_i'$.  So the map $\fy'$ satisfies the conditions
 ${\fy}' (x, 0) = Id$, ${\fy}' (x, 1) (P) = {\wP}'$ and for any $t\in[0,1]$,
${{\fy}' (x, t)}|_{{\wP}'} = {Id}_{{\wP}'}$.
\\
Define a map ${\fy}_i' : {P_i} \times [0, 1] \ra P_i$ by the
formula $${\fy}_i' (x, t) = S_i \circ {\fy}' ({S_i}^{-1} (x), t).$$
Each map ${\fy}_i'$ is a deformation retraction from $P_i$ to $S_i
({\wP}')$.
\\
Observe that the map ${\fy}_i'$ keeps all the vertices $S_i (A_k)$
of the polygon $P_i$ fixed. Therefore we can define a deformation
retraction ${\widetilde\fy} (x, t):{\wP} \times [0, 1] \to
 \bigcup\limits_{i = 1}^m{S_i ({\wP}')} ={\wP}$  by a formula
$${\widetilde\fy} (x, t) =  \fy_i' (x, t),
\mbox{\quad \rm  if  }x \in P_i$$
 The map $\widetilde\fy$ is well-defined and continuous
because if $P_i\bigcap P_j = \{S_i (A_k)\}=\{S_j(A_l)\}$ for some
$k$ and   $l$, then
${\fy}_i' (S_i (A_k), t) \equiv {\fy}_j' (S_j (A_l), t) \equiv S_i (A_k)$.
\\
 Moreover,${\widetilde\fy} (x, 0) =x$ on $\wP$, and
${\widetilde\fy} ({\wP}, 1) \equiv \bigcup\limits_{i = 1}^m S_i ({\wP}')$ and ${\widetilde\fy} (x, t)|_ {{\wP}''} \equiv Id$. \\ So ${\widetilde\fy} (x, t)$ is a deformation retraction from $\wP$ to ${\wP}''$.\\
Therefore, the set ${\wP}'' = \bigcup {S_i \circ S_j'(P)}$ is contractible.{\vse}

\begin{thm}\label{main}   Let $S$ be a polygonal tree system of similarities associated with $P$, and let $K$ be its attractor. Then $K$ is a dendrite.\end{thm}

{\dok} Let $T (A) =  {\bigcup}S_i (A)$ be the Hutchinson operator of the system $\eS$ and let ${\wP}^{(1)} = T(P), {\wP}^{(n + 1)} = T(\wP^{(n)})$.

By Lemma  \ref{trans}, each of the  sets  ${\wP}^{(n)}$ is a  contractible compact set, satisfying the inclusions  ${\wP}^{(1)} \NI  {\wP}^{(2)}\NI  {\wP}^{(3)} \ldots$. The diameter of connected components of the interior of   each ${\wP}^{(n)}$ does not exceed ${\diam}{P}\cdot {q^n}$, where $q = \max\Lip(S_i)$. Therefore the set $K = \bigcap {\wP}^{(n)}$ is contractible and has empty interior.
Since  the  system $\{P_i\}$ is  connected in the  sense of Definition \ref{conn}  by Kigami's theorem, the  attractor $K$ is connected, locally connected and  arcwise connected \cite[Theorem 1.6.2, Proposition 1.6.4]{Kig}.
Since any simple closed  curve in a contractible set $X$ on a plane bounds a disc in $X$ which has interior points, the set $K$ contains no simple closed  curve  and  therefore is a dendrite.
{\vse}

The dendrite $K$ lies in the polygon $P$, and its intersection with the sides of $P$ can be uncountable, 
or  even contain the  whole  sides of $P$. 
This is also true  for  any subpolygon $S_\bi(P)$.
 Nevertheless, all the  dendrite $K$ "squeezes" through the vertices  of each such subpolygon  $S_\bi(P)$, namely:

\begin{prop}\label{squeeze}
Let $\bj\in I^*$ be a multiindex. For any continuum $L\IN K$, whose intersection with both $ P_\bj$ and its exterior $\dot {CP_\bj}$ is nonempty, the set $\overline{L\mmm P_\bj}\cap P_\bj$ is a nonempty subset of the set $\{S_\bj(A_i), i=1,...,n\}$.

\end{prop} 
\dok Observe that for any polygon $P_\bj, \bj\in I^k$ the set  ${\wP}^{(k)}\mmm  \{S_\bj(A_i), i=1,...,n\}$ is not connected, and $P_\bj\mmm  \{S_\bj(A_i), i=1,...,n\}$ is its connected component, whose intersection with $K$ is equal to
$S_\bj(K\mmm \{A_i, i=1,...,n\})$. Therefore after deleting the vertices $\{S_\bj(A_i), i=1,...,n\}$, the continuum $L$ becomes disconnected too.\vse

\bigskip
\subsection{The main tree and ramification points}
Let ${\ga}_{ij}$ be the arc in $K$, connecting the vertices $A_i$ and $A_j$.

\begin{thm}\label{arcs}   The arcs ${\ga}_{ij}$ are the components of an invariant set of some multizipper $\eZ$.\end{thm}
{\dok} We say that the polygons $P_{i_1}, \ldots, P_{i_m}$ form a chain connecting $x$ and $y$, if $P_{i_1}\ni x,  P_{i_m} \ni y$ and $ P_{i_k}\bigcap  P_{i_l}$ is empty if $|l - k| > 1$ and is a common vertex of $ P_{i_k}$ and $ P_{i_l}$ when $| l - k| = 1$.
\\
For any $A_i, A_j$, there is a unique chain of polygons $P_{ijk}, k = 1, \ldots m_{ij}$ connecting  $A_i$  and $A_j$.
\\
Let $u (i, j, k)$ and $v (i, j, k)$ be such numbers that $S_{ijk}(A_u) = P_{ij{(k-1)}} \bigcap P_{ijk}= z_{ij{(k - 1)}}$ and  $S_{ijk}(A_v) =  P_{ij{k}} \bigcap P_{ij{(k+1)}}= z_{ijk}$, if $1 < k < m_{ij}$\\
$u (i, j,1) = A_i= z_{ij0}$ and $v (i, j, m_{ij})= A_j = z_{ijm_{ij}}$ \\
Then we have the following relations, $${\ga}_{ij} = \bigcup\limits_{i = 1}^{m_{ij}} S_{ijk}({\ga}_{u (i,j,k), v(i, j, k)}) = \bigcup\limits_{i = 1}^{m_{ij}}{\ga}_{ijk}.$$Therefore the system $\{S_{ijk}\}$ is a multizipper $\eZ$ with vertices  $z_{ijk}$.\\
Since each ${\ga}_{ijk}$ lies in $P_{ijk}$, \\
$$ {\ga}_{ijk} \bigcap {\ga}_{ijl} = \0,$$ if $|k - l| > 1$ and$${\ga}_{ijk} \bigcap {\ga}_{ijl} =\{ z_{ijk}\},$$ if $l = k \pm 1$.\\
Therefore, $\eZ$ satisfies the condition of Theorem \ref{Jor}.\\
So ${\ga}_{ij}$ are all Jordan arcs. {\vse}

\begin{dfn}\label{defmt}
The union $\hat\ga=\bigcup\limits_{i\neq j}\ga_{ij}$ is called   {\em  the main tree} of  the  dendrite $K$.
The ramification points of the  tree $\hat\ga$ are called {\em the main ramification points} of  the  dendrite $K$.
\end{dfn}

\noindent{\bf Example 2.3.}\bigskip\\
\begin{center}
\begin{tikzpicture}[line cap=round,line join=round,>=triangle 45,x=10.0cm,y=10.0cm]
\node[anchor=south west,inner sep=0] at (0,-.227) {
\includegraphics[width=.325\textwidth]{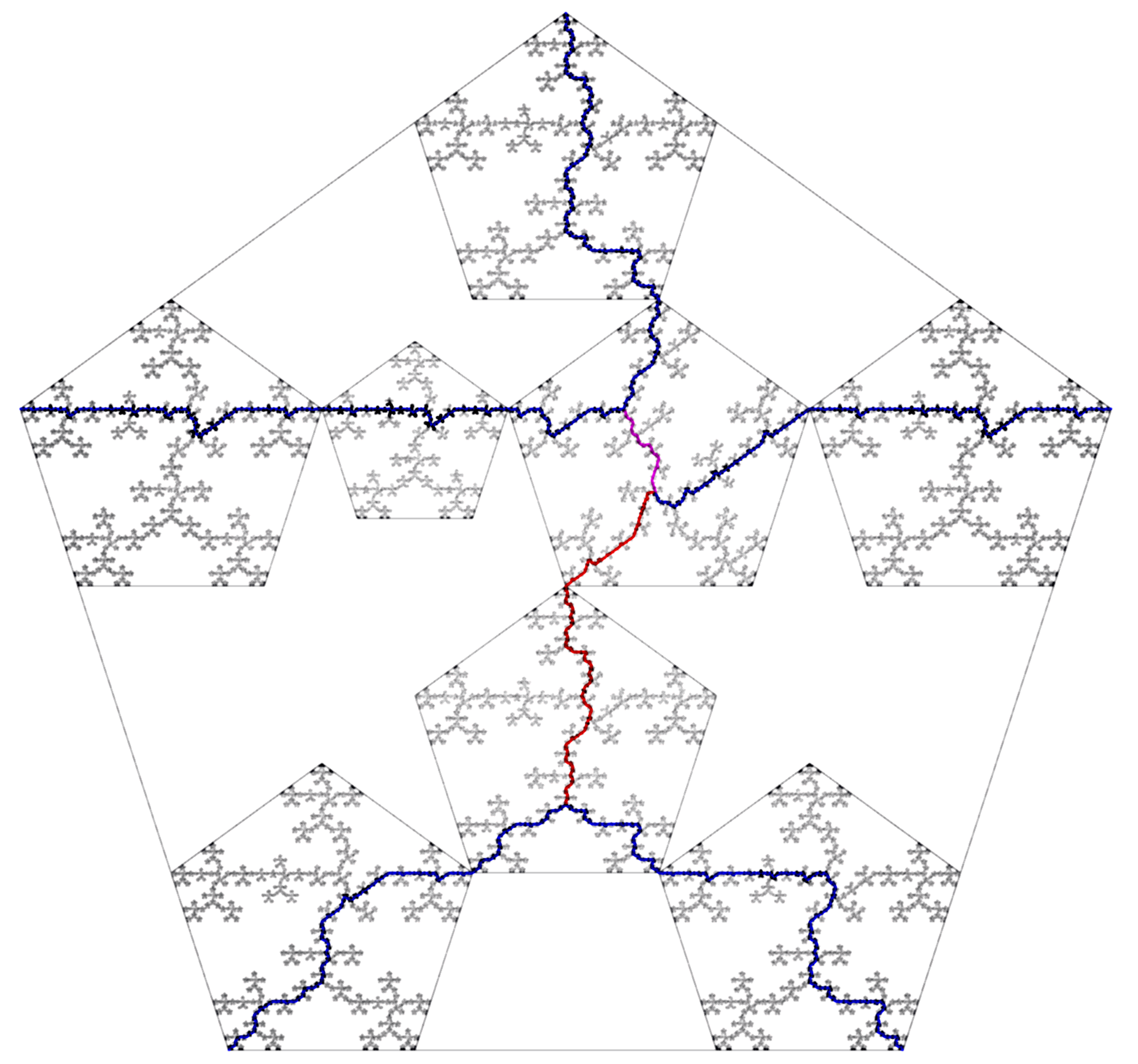}};
\node[anchor=south west,inner sep=0] at (.55,-.227) {\includegraphics[width=.54\textwidth]{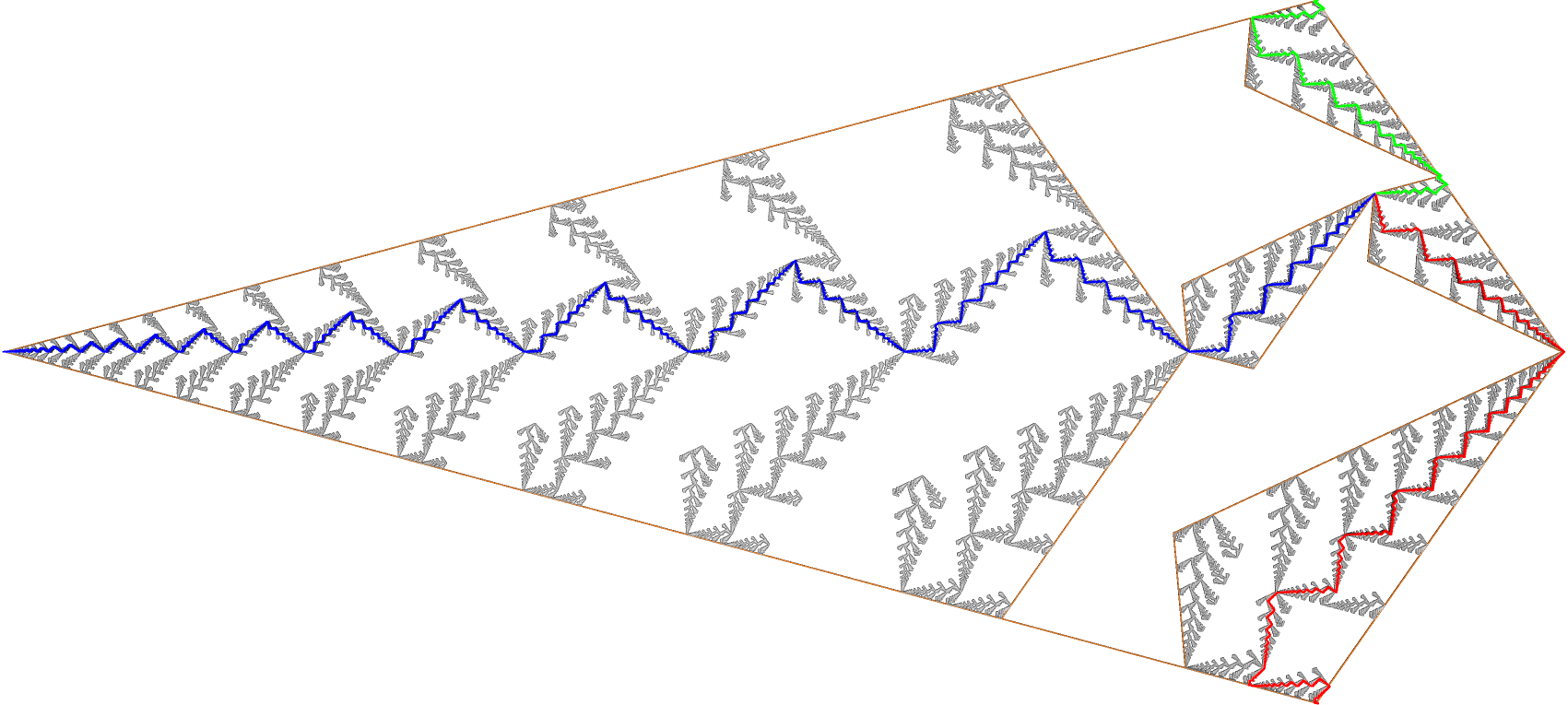}};
\draw [fill=red] (0.225,0.0135) circle (1.8pt);
\draw [fill=red] (0.235,-0.019) circle (1.8pt);
\draw [fill=red] (0.204,-0.133) circle (1.8pt);
\draw [fill=red] (1.15,-0.007) circle (1.8pt);
\end{tikzpicture}\\
{\small Two polygonal dendrites,  their main trees and main ramification points.}
\end{center}

{There is a simple way to know whether a point $x\in K$ lies in $\hat\ga$ and belongs to the set $CP(\hat\ga)$ of its cut points or to the set $EP(\hat\ga)$ of its end points:}\\

\begin{lem}\label{xintree} 
Let $x\in K$. (a) $x\in CP(\hat\ga)$ if and only if  there are  vertices $A_{i_1}$, $A_{i_2}$, not belonging to the same component of $K\mmm\{x\}$; (b) $x\in EP(\hat\ga)$ iff $x$ is a vertex and $x\notin CP(\hat\ga)$.
\end{lem}
\dok  First part of (a) is obvious. Since  the union $\ga_{xA_{i_1}}\cup\ga_{xA_{i_2}}$ is a Jordan arc, it is equal to $\ga_{i_1i_2}$. So $x$ is a cut point of $\ga_{i_1i_2}$, and therefore of $\hat\ga$. 
To check (b), suppose $x\in\hat\ga$ is not a vertex, then $x$ lies in some $\ga_{i_1i_2}$, so it is  a cut point of $\hat\ga$. The second part of (b) is obvious.
\vse

There are   points in $K$ for which their order in $K$ and in $\hat\ga$ is the same:

\begin{lem}\label{eqord} 
Let $x\in CP(K)$. If each component $C_l$ of $K\mmm\{x\}$ contains a vertex of $P$, then $Ord(x,K)$ is finite and $Ord(x,K)=Ord(x,\hat\ga)$\end{lem}
\dok The number of components of $K\mmm\{x\}$ is not greater than $n$, so it's finite. 
Let  $ C_l,l=1,...,k, k=Ord(x,K)$ be the components of $K\mmm\{x\}$.
 By Lemma  \ref{xintree}, $x\in\hat\ga$. It also follows from Lemma \ref{xintree} that two vertices $A_{i_1}$ and $A_{i_2}$ lie in the same component $C_l$ if and only if $x\notin\ga_{i_1i_2}$. Therefore, all the vertices of $P$, belonging to the same component $C_l$ of $K\mmm \{x\}$, belong to the same component of $\hat\ga\mmm \{x\}$. Therefore $Ord(x,{\hat\ga})=Ord(x,K)$.\bigskip\vse

For $x\in\rr$, we denote by $\lceil x\rceil$  {\em the ceiling of $x$}, or the minimal integer $n$ which is  greater or equal to $x$.

\begin{prop}\label{comp} 
a) For any $x\in\hat\ga$, $\hat\ga=\bigcup\limits_{ j=1}^n\ga_{A_jx}$.\\
b) $A_i$ is a cut point of $\hat\ga$, if there are $j_1,j_2$ such that $\ga_{j_1i}\cap\ga_{j_2i}=\{A_i\}$;\\
c) the only end points of $\hat\ga$ are the vertices $A_j$ such that $A_j\notin CP(\hat\ga)$;\\
d) if $\#\pi^{-1}(A_i)=1$, then $Ord(A_i,\hat\ga)=Ord(A_i,K)\le n-1$,
otherwise $Ord(A_i,K)\le (n-1)(\left\lceil{\dfrac{\te_{max}}{\te_{min}}}\right\rceil-1)$, where $\te_{max},\te_{min}$ are maximal and minimal values of vertex angles of $P$.
 \end{prop}
\dok   For any $j_1,j_2$, $\ga_{j_1j_2}\IN \ga_{A_{j_1}x}\cup\ga_{A_{j_2}x}$, which implies a). Repeating argument of Lemma \ref{xintree}, we see that $A_i$ is a cut point of $\ga_{i_1i_2}$ and therefore of $\hat\ga$, thus proving b). 
If $x\in\hat\ga$ is not  a vertex, then for some $j_1,j_2$, $x\in\ga_{j_1j_2}$, so $x$ is a cut point of $\ga_{j_1j_2}$ and therefore of $\hat\ga$, which implies c).

Let  $\{C_l,l=1,...,k\}$ be some set of components of $K\mmm\{ A_i\}$. Since $\{A_i\}$ is the intersection of unique nested sequence of polygons $P_{j_1}\NI P_{j_1j_2}\NI...\NI P_{j_1..j_s}..$, 
there is such $s$, that $\diam P_{j_1..j_s}<\diam C_i$ for any $i=1,...,k$. Then, by Proposition \ref{squeeze}, each $C_l$ contains some vertex of $P_{j_1..j_s}$, different from $A_i$, therefore $k\le n-1$ so $Ord(A_i,K)\le n-1$ is finite. So we can suppose that we took $k=Ord(A_i,K)$ initially and $\{C_1,...,C_k\}$ was the set of all components of $K\setminus \{A_i\}$.  

Let $\bj=j_1..j_s$ and $A_i=S_\bj(A')$.
The sets $C_l\cap P_\bj$ are the components of $K_\bj\mmm\{A_i\}$. Since $(K\cap P_\bj)\mmm\{A_i\}=S_\bj(K\mmm\{A'\})$, there are  $k$ components $C_l'$ of $K\mmm \{A'\}$, such that $S_\bj(C_l')=C_l\cap P_\bj$. Since  each set $C_l'$ contains the vertices of $P$, by Lemma \ref{eqord}, $Ord(A',{\hat\ga})=Ord(A',K)=Ord(A_i,K)\le n-1$.

Suppose $\#\pi^{-1}(A_i)>1$, and let $P_{j_1}\NI P_{j_1j_2}\NI..$ and $P_{j_1'}\NI P_{j_1'j_2'}\NI...\NI P_{j_1'..j_s'}..$ be two different nested sequences of polygons whose intersection is $A_i$. 
For any two polygons $P_\bj, P_{\bj'}$ either their intersection is $A_i$  or one of these polygons contains the other. Therefore, there is some $k$ such that $P_{j_1..j_s}=P_{j_1'..j_s'}$ for $s<k$ and $P_{j_1..j_s}\cap P_{j_1'..j_s'}=\{A_i\}$ for $s\ge k$.
Since the  vertex angles of respective polygons at $A_i$ form  a decreasing sequence assuming finite set of values, 
both sequences of these values are eventually constant. 
These final values are greater or equal to $\te_{min}$. 
Therefore, there is a finite number of polygons $P_{\bj^k}\ni A_i$, 
whose pairwise intersections are $\{A_i\}$, such that any  other  polygon $P_{\bj'}$, containing $A_i$, either contains one of them,
 or is contained in some $P_{\bj^k}$ and has the same vertex angle at $A_i$.
  Then $Ord(A_i,K)=\sum Ord(A_i,P_{\bj^k})=\sum Ord(A_i,S_{\bj^k}(\hat\ga))$. The number of polygons $P_{\bj^k}$  is not greater than $\left\lceil{\dfrac{\te_{max}}{\te_{min}}}\right\rceil-1$, therefore $Ord(A_i,K)\le (n-1) \left\lceil{\dfrac{\te_{max}}{\te_{min}}-1}\right\rceil$
\bigskip\vse

\noindent{\bf Example 2.4.}\bigskip\\
\begin{tikzpicture}[line cap=round,line join=round,>=triangle 45,x=10.0cm,y=10.0cm]
\node[anchor=south west,inner sep=0] at (0,-.227) {
   \includegraphics[width=.79\textwidth]{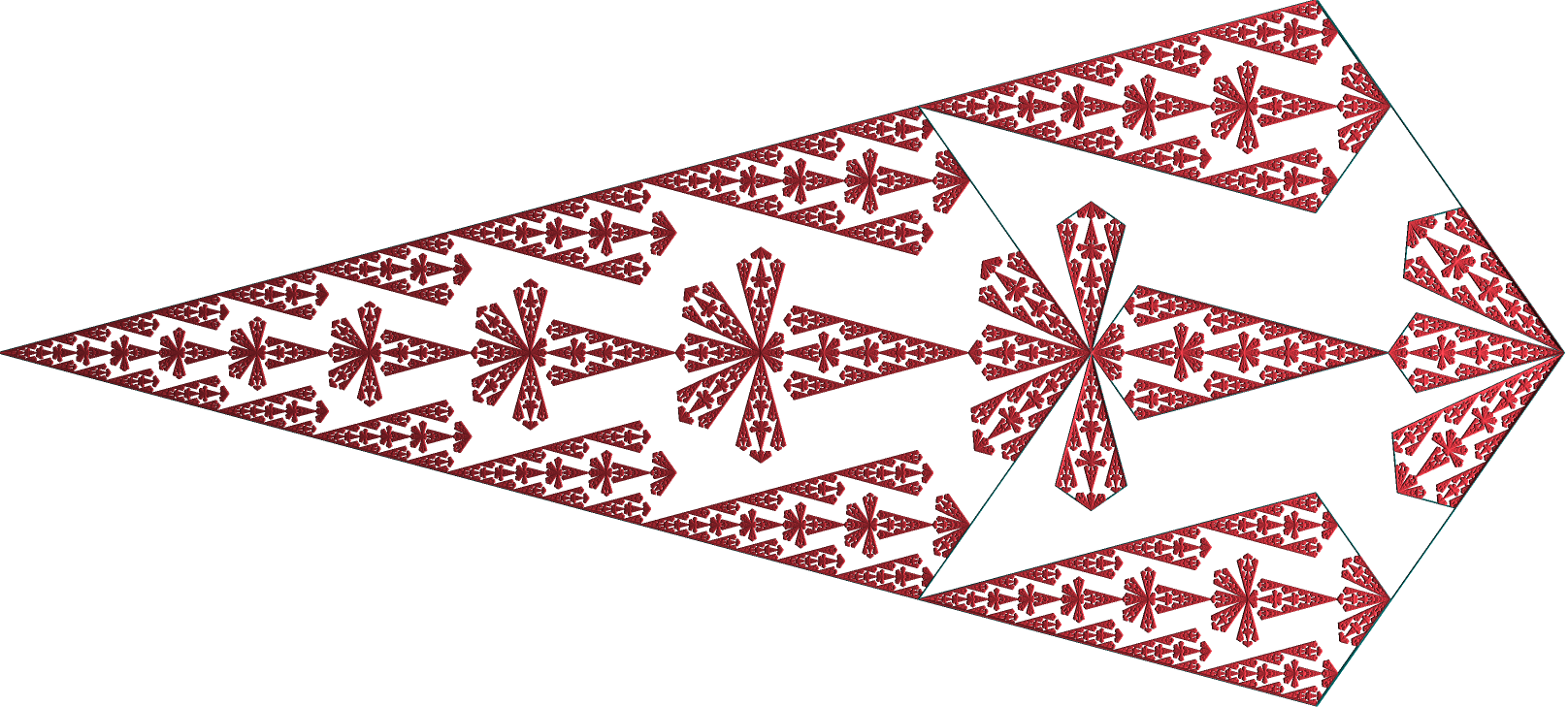}};
%\draw (0,0.15) node[anchor=south west] {\bf $Ord(x,K)$ for some points:};
\draw (0,0.05) node[anchor=north east] {\bf A};
\draw (1,0.05) node[anchor=north west] {\bf C};
\draw  (0.842,0.22) node[anchor=south west] {\bf D};
\draw  (0.842,-0.22) node[anchor=north west] {\bf B};
\draw (.33,-.2) node [anchor=south east] {\bf E};
\draw[->,color=blue!50!black,ultra thick](.3,-.2)--(.698,0);
\draw[color=green!50!black,ultra thick] (0.,0.) -- (0.842,0.226) -- (1.,0.) -- (0.842,-0.226) -- (0.,0.) -- (0.588,0.157) -- (0.698,0.) -- (0.588,-0.157) -- (0,0);
\draw[color=green!50!black] (0.588,0.157) -- (0.842,0.226) -- (0.89,0.157) -- (0.842,0.089) -- cycle;
\draw[color=green!50!black] (0.588,-0.157) -- (0.842,-0.0893) -- (0.89,-0.157) -- (0.842,-0.226) -- cycle;
\draw[color=green!50!black] (1.,0.) -- (0.8975,0.048) -- (0.901,0.084) -- (0.935 ,0.093) --  (1.,0.) -- (0.905,-0.025) -- (0.888,0.) -- (0.905,0.025) --  (1.,0.) -- (0.93,-0.099) -- (0.893,-0.089) -- (0.89,-0.051) -- (1,0);
\draw[color=green!50!black] (0.888,0.) -- (0.728,-0.043) -- (0.698,0.) -- (0.728,0.043) -- cycle;
\draw[color=green!50!black] (0.698,0.) -- (0.676 ,0.081 ) -- (0.698,0.096 ) -- (0.72,0.081) --  (0.698,0.) -- (0.72 ,-0.085) -- (0.698,-0.101) -- (0.675,-0.085) -- cycle;
\end{tikzpicture}\\
{\small  A polygonal system, generated by 9 maps of a quadrilateral ABCD with vertex angles  30, 110, 110 and 110 degrees. The main tree $\hat\ga$ is the union of line segments AB, AC and AD. For the vertex A, $Ord(A,K)=Ord(A,\hat\ga)=3$.
The vertices B and D have order 1 both in $\hat\ga$ and in the dendrite K. The vertex C has $Ord(C,\hat\ga)=1$ but $Ord(C,K)=9$.
The point E has the maximal order 24. By Theorem \ref{order}, for this type of polygon, maximal possible order may be 33. }

\begin{thm}\label{order}  For each cut point $y\in K$ there is $S_\bi$  such that for some $x\in\hat\ga$,
$y=S_\bi(x)$.   If $x$ is not a vertex of $P$,   $Ord(y,K)=Ord(x,{\hat\ga})$. Otherwise, there are 
multiindices $\bi_k, k=1,..,s$ and vertices $x_1,...,x_s$, 
such that for any  $k$, $S_{\bi_k}(x_k)=y$, for any $l\neq k$,  $S_{\bi_k}(P)\cap S_{\bi_l}(P)=\{y\}$ and $ Ord(y,K)=\sum\limits_{k=1}^s Ord(x_k,{\hat\ga})\le(n-1)\left(\left\lceil{\dfrac{2\pi}{\te_{min}}}\right\rceil-1\right)$.
\end{thm}

{\bf Proof.} Let  $\{C_1,...,C_k\}, k>1$, be some set of the components of $K\setminus \{y\}$. 
Take $0<\rho<\min\limits_{i=1,...,k}\diam(C_i)$. 
Let $\bj\in I^*$ be a  multiindex such that $P_\bj\ni y$ and $\diam(P_\bj)\le\rho$ and let $y=S_\bj(x)$.

Suppose the point $x$ is not a vertex of the polygon $P$. 
Then $y\in\dot P_\bj$ and the sets $C_i\cap P_\bj$ are the components of $K_\bj\mmm\{y\}$. Since $(K\cap P_\bj)\mmm\{y\}=S_\bj(K\mmm\{x\})$, there are    $k$ components $C_i'$ of $K\mmm \{x\}$, such that $S_\bj(C_i')=C_i\cap P_\bj$. By Proposition \ref{squeeze},  each set $C_i'$ contains the vertices of $P$, therefore $k\le n$ and $Ord(y,K)\le n$. So we can suppose that we took $k=Ord(y,K)$ initially and $\{C_1,...,C_k\}$ was the set of all components of $K\setminus \{y\}$. 
Since each set $C_i'$ contains the vertices of $P$, by Lemma \ref{eqord}, $Ord(x,{\hat\ga})=Ord(x,K)=Ord(y,K)$.

The proof of the last part repeats the proof of  d) in Proposition \ref{comp}.\vse

\begin{cor}\label{pcf}
Let $(P,\eS)$ be a polygonal tree system and $K$ be its attractor.
(i) For any $x\in K$, the set $\pi^{-1}(x)$ contains no more elements than $(n-1)\left(\left\lceil{\dfrac{2\pi}{\te_{min}}}\right\rceil-1\right)$;\\
(ii)  The system $\eS$ is post-critically finite.
\end{cor}
\dok (i) was proved in previous Theorem. Since post-critical set is contained in $\pi^{-1}(\{A_1,...,A_n\})$, it is finite.\vse

\bigskip

\subsection{Metric properties of polygonal dendrites.}

  Following \cite{TV}, we remind that for $c\ge 1$, a set $A \IN R^n$ is of {\em c-bounded turning} if
each pair of points $a, b\in A$ can be joined by a continuum $F\IN A$ with diameter
$\diam(F)\le c|a-b|$. In this subsection we prove that a dendrite $K$, defined by a polygonal tree system, is of c-bounded turning for some $c\ge 1$.

 \begin{lem}\label{rho} Let $\{P, \eS\}$ be a polygonal tree system.  There is such $\rho$ that\\ (i) for any vertex $A$,   $V_\rho(A) \bigcap P_k \ne \0 \Rightarrow P_k \ni A$;\\ (ii) for any $x, y \in P$ such that there are $P_k, P_l: x \in P_k, y \in P_l$ and $P_k \bigcap P_l
   = \0, d (x, y) \ge \rho$. \vse  \end{lem}
   Let $\al$ denote the minimal angle between the sides of polygons $P_i, P_j$, having common vertex.

 \begin{lem}\label{ga_xA} For any vertex $A$ of $P$ and for any $x \in K \setminus \{A\}$, $$\dfrac{\diam \ga_{Ax}} {d (x, A)} \le \dfrac {\diam P}
   {\rho}$$\end{lem}
   \dok There are such $i_1,  \ldots, i_{k+1}$ that $A \in S_{i_1\ldots i_{k+1}} (P)$ and\\ $x \in S_{i_1\ldots i_k} (P) \setminus
   S_{i_1\ldots i_{k+1}} (P)$. Let $x' = S^{-1}_{i_1\ldots i_k} (x)$  and $A' =   S^{-1}_{i_1\ldots i_k} (A)$. Then $x' \in P \setminus
   P_{i_{k+1}}$ and $A' \in P_{i_{k+1}}$,  so $d (x', A') \ge \rho$, and $\dfrac{\diam \ga_{x' A'}} {d(x', A')} \le \dfrac {\diam P} {\rho}$. Since $S_{i_1\ldots i_k}(\ga_{x' A'})=\ga_{xA}$, we get $\dfrac{\diam \ga_{x A}} {d(x, A)} \le \dfrac {\diam P} {\rho}.$ \vse

   \begin{lem}\label{ga_xy} If $x \in S_k(K), y \in S_l(K), P_k\cap P_l = A$ and $x \ne y$, then $$\dfrac{\diam \ga_{xy}}{d(x, y)} \le
   \dfrac{\diam P}{\rho
  \sin{(\al/2)}}$$.\end{lem}
\dok
$\dfrac {d(x, y)} {d (x, A) + d (A, y)}\ge \dfrac {\sqrt {{d (x, A)}^2 + {d (A, y)}^2-2 d (x, A)d (A, y)\cos{\alpha}}}  {d (x, A) + d (A, y)}$.\\

The minimum value for the right side of equation over all $d(x,A),d(y,A)$ is $\sin\al/2$, while, by Lemma  \ref{ga_xA},
 $$\dfrac {d (x, A) + d (A, y)}{\diam \ga_{xy}} \ge \dfrac {\rho}{\diam P}\eqno{(2)}$$\\
Therefore we have
$\dfrac{\diam \ga_{xy}}{d(x, y)} \le
   \dfrac{\diam P}{\rho
  \sin{(\al/2)}}$. \vse

   \begin{lem}  For any $x, y \in K$,
   $\dfrac{\diam \ga_{xy}}{d(x, y)} \le  \dfrac {\diam P}{\rho
   \sin{(\al/2)}}$.  \end{lem}
   \dok
   There are such $i_1, \ldots, i_k, i_{k+1}$ that $x \in S_{i_1 \ldots i_{k+1}} (P)$ and $y \in   S_{i_1 \ldots i_k}{ (P \setminus P_{i_{k
   +1}})}$. Let $x' = S^{-1}_{i_1\ldots i_k} (x), y' = S^{-1}_{i_1\ldots i_k} (y)$. Suppose $y' \in P_l$.\\

   If $P_l \bigcap P_{i_{k+1}} =
   \0$, then $\dfrac{\diam \ga_{x'y'}} {d(x', y')} \le \dfrac{\diam P} { \rho}$.\\

    If $P_l$ and $P_{i_{k+1}}$ have a common vertex, then $\dfrac{\diam \ga_{x' y'}} {d(x', y')} \le \dfrac{\diam P}{ \rho\sin\al/2}$.\\

     Thus we have,
  $$\dfrac{\diam \ga_{xy}}{d(x, y)} \le  \dfrac {\diam P}{\rho
   \sin\al/2}.\quad \blacksquare$$
 From previous  three Lemmas  we immediately get the following
 \begin{thm}\label{c-BT} The attractor $K$  of a   polygonal tree  system $\eS$   is a continuum with bounded turning.\ \vse
 \end{thm}

\subsection{Morphisms of polygonal dendrites}

In the following Theorem we admit that the enumeration of the vertices of the polygons $P$ and $P'$ needs not follow any order, and all permutations of indices are allowed.

  \begin{thm}\label{morphism} Let dendrites $K,  K'$ be the  attractors of polygonal tree systems $\eS = \{S_1,  S_2,  \ldots S_m\}$ and ${\eS'}=\{S'_1,  S'_2,  \ldots S'_m\}$ associated with polygons $P,  P'$ whose vertices $A_1, ...,A_n$ and $A_1',...,A_n'$ satisfy the conditions\\

(i) For any $i,j=1,..., n$, $S_k (A_i) = A_j$ iff $S'_k (A'_i)= A'_j$;\\

(ii) For any $i,j=1,..., n$ $S_{k_1} (A_{i}) = S_{k_2} (A_j)$ iff $ S'_{k_1} (A'_i) =
S'_{k_2} (A'_j)$. \\

Then there is a  bi-H\"older  homeomorphism $\psi: K \to K'$ such that for any $i=1,...,m$, $  \psi \circ S_i =
S'_i \circ \psi$.\end{thm}

\dok

{\bf  1.}  {\em The condition (i)  implies that  for any multiindex $\bk=k_1 k_2... k_l\in I^*$ the equality
$S_\bk(A_i)=A_j$ holds iff $S'_\bk(A'_i)=A'_j$}. \\

Indeed, it's true for $l=1$; proceeding by induction, let the condition (i) be true for any $k_1  k_2... k_l \in I^l$ and $i, j \in \{1, \ldots, n\}$, i.e.
$$S_{k_1 \ldots k_l} (A_i) = A_j \Longleftrightarrow S'_{k_1 \ldots k_l} (A'_i) = A'_j$$

Suppose for some $k_1 k_2...k_{l+1} \in I^{l+1}$ and some vertices $A_i, A_j$ we have  $S_{k_1 k_2 \ldots k_{l+1}} (A_i) = A_j.$\\

Consider the point $S_{k_2 \ldots k_l k_{l+1}} (A_i)=S_{k_1}^{-1}(A_j)$.
This point is some vertex $A_{i_1}$ of $P$.
Since the multiindex  $k_2, \ldots, k_l, k_{l+1}$ is of length $l$,
 $S'_{k_2 \ldots k_l k_{l+1}} (A'_i)=A'_{i_1}$ by induction hypothesis.
At the same time,  $S'_{k_1} (A'_{i_1}) = A'_j$ .
Therefore  $S'_{k_1 k_2 \ldots k_l k_{l+1}} (A'_i) = A'_j$. \\

{\noindent\scriptsize \begin{tikzpicture}[line cap=round,line join=round,>=triangle 45,x=6.0cm,y=6.0cm]
\draw (0.,0.) -- (0.842,0.2256) -- (1.,0.) -- (0.842,-0.2256) -- cycle;

%polygons and their filling
\filldraw[color=black,fill=blue,fill opacity=0.2] (0.,0.) -- (0.5377,0.144) -- (0.638,0.) -- (0.5377,-0.144) -- cycle;
\filldraw[color=black,fill=blue,fill opacity=0.2] (0.638,0.) -- (0.773,0.0359) -- (0.798,0.) -- (0.773,-0.035) -- cycle;
\filldraw[color=black,fill=green,fill opacity=0.2] (0.773,0.035) -- (0.788,0.21) -- (0.8420,0.2256) -- (0.873,0.18) -- cycle;
\filldraw[color=black,fill=blue,fill opacity=0.2] (0.798,0.) -- (0.9681,0.0452) -- (1.,0.) -- (0.968,-0.0455) -- cycle;
\filldraw[color=black,fill=red,fill opacity=0.1](0.773,-0.0359) -- (0.873,-0.18) -- (0.842,-0.2256) -- (0.788,-0.21) -- cycle;

%vertices A1-A4$
\draw (0,-0.03) node {$ A_1$} (1.02,-0.03) node {$A_3$}
(0.83,0.27) node {$ A_4$} (0.83,-0.27) node {$ A_2$};

\draw (0.07,0) node { $ 1$} [red!60!black](0.52,-0.12) node {$   2$}  [blue!60!black](0.606,0) node{$  \ 3$} (0.52,0.12)[green!60!black] node {$  4$};

\draw (0.76,-0.06) node { $\bf 1$} [red!60!black](0.77,-0.19) node {$  \bf 2$} [blue!60!black](0.83,-0.2) node{  $\bf 3$} (0.87,-0.15)[green!60!black] node {$  \bf 4$};

\draw (0.76,0.06) node { $\bf 1$} [red!60!black](0.87,0.15) node {$  \bf 2$} [blue!60!black](0.83,0.2) node{  $\bf 3$}(0.77,0.19) [green!60!black]node {$  \bf 4$};

%maps S1-S5
\draw (0.3,0) node {$  S_1$}(0.74,0)node {$  S_5$}(0.92,0)node{$S_3$}(0.81 ,-0.15 ) node{$ S_2$}(0.81 ,0.15 ) node{$ S_4$};

\end{tikzpicture}  \begin{tikzpicture}[line cap=round,line join=round,>=triangle 45,x=6.0cm,y=6.0cm, baseline=-12mm ]

\filldraw[color=black,fill=blue,fill opacity=0.005] (0.,0.) -- (1.,0.) -- (0.71,0.29) -- (0.29,0.29) -- cycle;
\filldraw[color=black,fill=blue,fill opacity=0.2] (0.,0.) -- (0.333,0.) -- (0.23,0.0976) -- (0.097,0.0976) -- cycle;
\filldraw[color=black,fill=blue,fill opacity=0.2] (0.333,0.) -- (0.666,0.) -- (0.569,0.0976) -- (0.4309,0.0976) -- cycle;
\filldraw[color=black,fill=blue,fill opacity=0.2] (0.666,0.) -- (1.,0.) -- (0.902,0.0976) -- (0.764,0.0976) -- cycle;
\filldraw[color=black,fill=red,fill opacity=0.1] (0.569,0.097) -- (0.707,0.292) -- (0.609,0.276) -- (0.552,0.195) -- cycle;
\filldraw[color=black,fill=green,fill opacity=0.1] (0.43096,0.0976) -- (0.292,0.2928) -- (0.276,0.1952) -- (0.333,0.114) -- cycle;

\draw (0.06,0.025) node{$ \bf1$}[blue!60!black](0.27,0.025) node {$ \bf 3$} [red!60!black](0.22,0.075) node{$ \bf 2$} (0.105,0.075) node [green!60!black]{$ \bf 4$};

\draw (0.393,0.025) node {$ \bf 1$}[blue!60!black](0.61,0.025) node{$ \bf 3$}[red!60!black](0.56,0.075) node {$ \bf 2$} (0.438,0.075) node [green!60!black]{$ \bf 4$};

\draw (0.727,0.025) node{$ \bf 1$} [blue!60!black](0.95,0.0255) node{$ \bf 3$} [red!60!black](0.89,0.075) node{$ \bf 2$} (0.772,0.075) node[green!60!black] {$ \bf 4$};

\draw (0.58,0.13) node {$ \bf 1$} [blue!60!black](0.67,0.265) node {$\bf 3$} [red!60!black](0.62,0.26) node {$ \bf 2$} (0.53,0.2) node[green!60!black] {$ \bf 4$};

\draw (0.4,0.12) node {$ \bf 1$} [blue!60!black](0.3,0.255) node {$ \bf 3$} [red!60!black](0.255,0.2) node  {$ \bf 2$} (0.35,0.13) node[green!60!black] {$ \bf 4$};

\draw (0.16,0.05) node {$ S_1$}(0.61,0.2) node{$ S_2$} (0.84,0.05) node {$ S_3$} (0.33,0.2) node {$ S_4 $} (0.51,0.05) node{$S_5$};

\draw (0,-0.03) node {$ A_1$} (1.,-0.03) node {$ A_ 3$} (0.69,0.33) node{$ A_2$} (0.27,0.33) node {$A_4$};

\end{tikzpicture}}

\begin{center} {\small Permutation of the vertices   defining an isomorphism of two polygonal tree systems.  The respective attractors are shown below.}
\end{center}

\noindent\includegraphics[width=.45\textwidth]{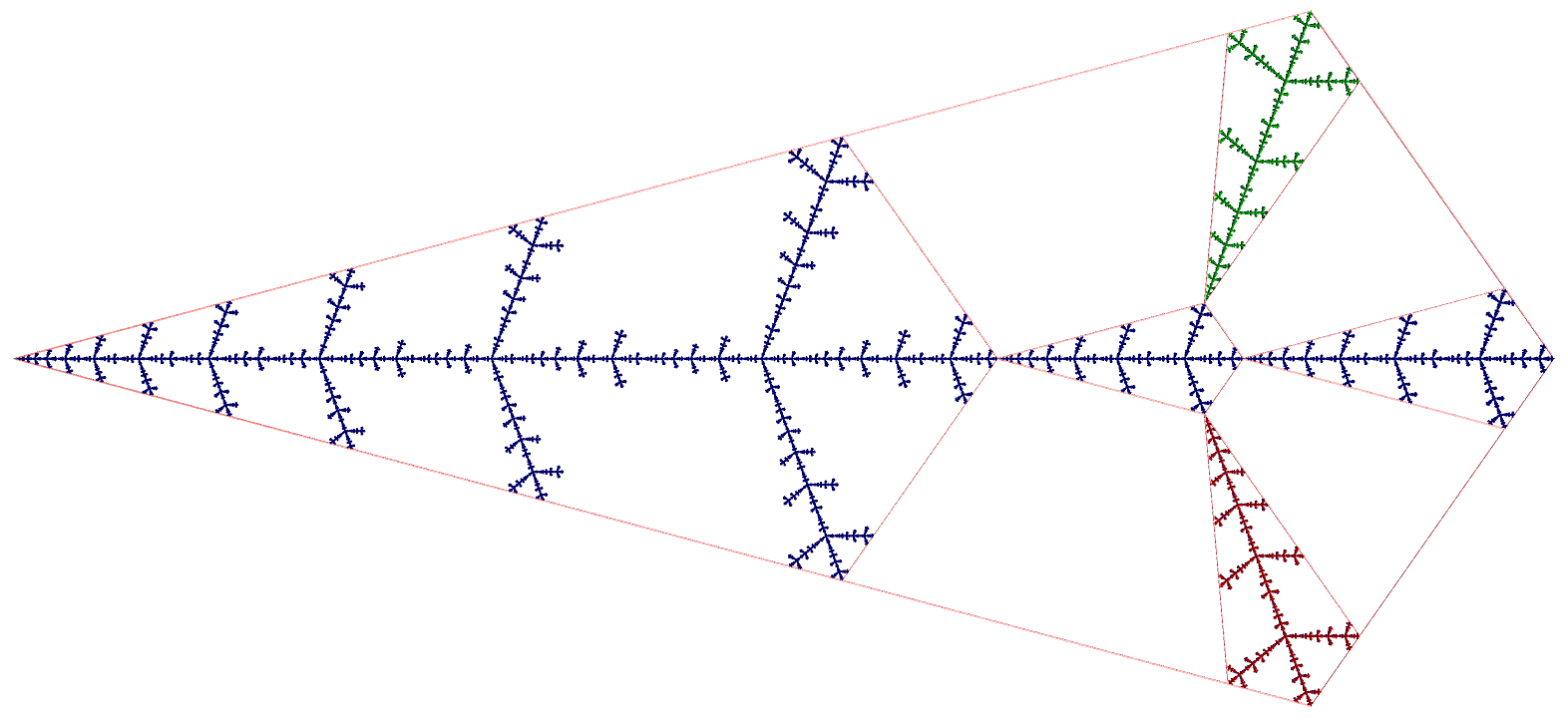} \qquad \includegraphics[width=.45\textwidth]{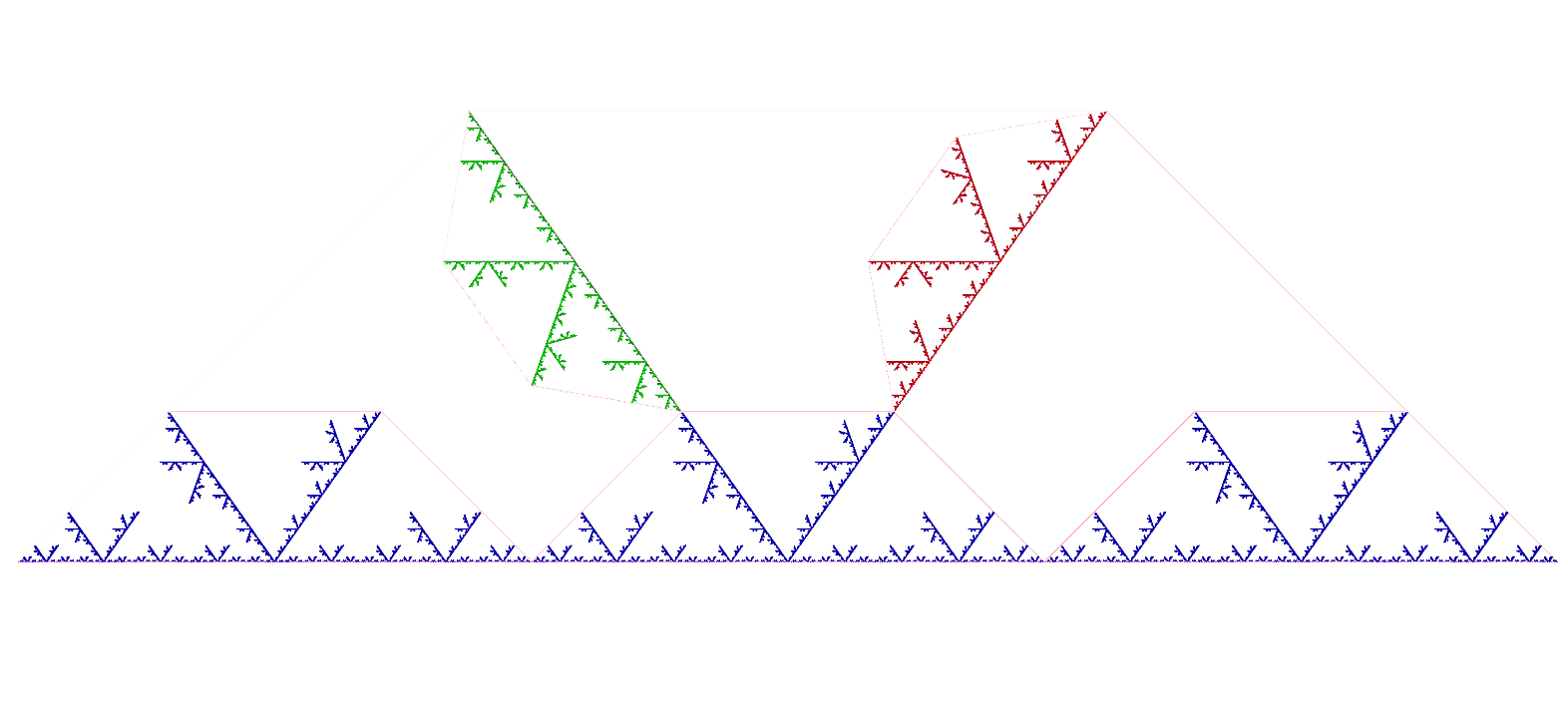}

\vspace{.3 cm}

{\bf 2.}  {\em The condition (ii)  implies that  for any multiindices $p_1...p_k$ and $q_1...q_l$ the equality
$S_{p_1 \ldots p_k }  (A_i) = S_{q_1 \ldots q_l} (A_j)$ holds iff $S'_{p_1 \ldots p_k }  (A'_i) = S'_{q_1 \ldots q_l} (A'_j)$}. \\

Suppose for some multiindices $p_1...p_k$ and $q_1...q_l$  and   vertices $A_i,  A_j$, \\ $S_{p_1 \ldots p_k }  (A_i) = S_{q_1 \ldots q_l} (A_j)$.\\
Rewrite it as  $S_{p_1} (S_{p_2 \ldots p_k }  (A_i) )= S_{q_1} (S_{q_2 \ldots q_l} (A_j))$.\\
Since $S_{p_2 \ldots p_k} (A_i)=S_{p_1}^{-1}(A_j)$, this point  must be   some vertex $A_{i_1}$ of $P$.
 Similarly, we also have $S_{q_2 \ldots q_l} (A_j) = A_{j_1}$.\\

 From (i) it follows that $S'_{p_2 \ldots p_l} (A'_i) = A'_{i_1}$  and  $S'_{q_2 \ldots q_l} (A'_j) = A'_{j_1}$ and
 from  $S_{p_1} (A_{i_1}) = S_{q_1} (A_{j_1})$  by (ii) it follows that  $S'_{p_1} (A'_{i_1}) = S'_{q_1} (A'_{j_1})$.

 Therefore, we have  $S'_{p_1 \ldots p_k }  (A'_i) = S'_{q_1 \ldots q_l} (A'_j)$.\\

{\bf 3.} {\em  There is a bijection $\fy: K\to K'$, such that for any $i\in I$, $\fy\cdot S_i=S'_i\cdot\fy$}.\\

Consider the index maps $\pi:\io\to K$ and $\pi':\io\to K'$.

Suppose for some $\Ep=p_1p_2p_3....\in \io$ and $\Eq=q_1q_2q_3....\in \io$,
$\pi(\Ep)=\pi(\Eq)=\{x\},  x\in K$.

Then for  any $k,l\in \nn$,  $P_{p_1...p_k}\cap P_{q_1...q_l}=\{x\}$.

Therefore, for any $k,l$ there are  such vertices
$A_{i_k}, A_{j_l}$ that $S_{p_1...p_k}(A_{i_k})= S_{q_1...q_l}(A_{j_l})=x$.
 Then, for any $k,l$, $S'_{p_1...p_k}(A'_{i_k})= S'_{q_1...q_l}(A'_{j_l})$.
 These equations imply  the points $S'_{p_1...p_k}(A'_{i_k})$ and $S'_{q_1...q_l}(A'_{j_l})$ coincide for all  $k,l$ and therefore
$\bigcap\limits_{k=1}^\8 P'_{p_1...p_k}=\bigcap\limits_{l=1}^\8 P'_{q_1...q_l}$.
 Applying this to all possible sequences $\Ep\in\pi^{-1}(x)$,
we obtain that $\pi'(\pi^{-1}(x))$ is a unique point, which we denote as $x'$.

Denote the map $\pi'\cdot\pi^{-1}:K\to K'$ by $\fy$.
Since the same argument shows that $\pi\cdot{\pi'}^{-1}:K'\to K$ is the inverse map to $\fy$,  the map $\fy $ is a bijection.

Since $\pi$ and $\pi'$ are compatible  with  the self-similar  structure on $\io, K$ and $K'$, the same is true for $\fy=\pi'\cdot \pi^{-1}$.

   \bigskip

 {\bf 4.} {\em The maps $\fy$ and $\fy^{-1}$ are H\"older continuous.} \\

   Denote $r_i = \Lip S_i, r_i' = \Lip S_i', \beta = \min \limits _{i = 1, \ldots, m} {\dfrac {\log {r'_i}} {\log {r_i}}}$, $\beta' = \min \limits _{i = 1, \ldots, m} {\dfrac {\log {r_i}} {\log {r'_i}}}$. Let also $|P| , |P'| $ be the diameters of  $P$  and $P'$  respectively. Let $\rho$ and $\rho'$ denote  the minimal distances specified by Lemma  \ref{rho}  for  the systems $\eS$ and $\eS'$ respectively and let $\al$, $\al'$ be respective minimal angles.

   Observe that for any multiindex ${\bf i} = i_1, \ldots, i_k, r_{\bf i}' \le r_{\bf i}^{\beta}$\\

  Take some $x, y \in K$. There is a multiindex  $i_1 \ldots i_k$ such that  $\{x,y\}\IN P_{i_1 \ldots i_k}$ and  for any $i_{k+1}$,  $\{x,y\}\not\subset P_{i_1 \ldots i_k i_{k+1}}$.
  Then there are two possibilities:\\

   a) For some pair of multiindices, $i_1 \ldots i_k j$ and   $i_1\ldots i_k  l$,  \\
$ P_{i_1 \ldots i_k j}\cap
   P_{i_1\ldots i_k  l}=\0$, \quad
 $x \in P_{i_1 \ldots i_k j}$ \quad  and  $y \in
   P_{i_1\ldots i_k  l}$.  \\

Then $d (x,y) \le r_{i_1 \ldots i_k} |P|$,
    while by Lemma  \ref{rho},  $d (x, y) \ge r_{i_1 \ldots i_k} \rho$.\\ In this case, $ r_{i_1 \ldots i_k} \rho < d (x,y) \le r_{i_1 \ldots i_k}|P|$. \\
   The same way,  for the  system $\eS'$  we have   $ r'_{i_1 \ldots i_k} \rho'_1 < d (x',y') \le r'_{i_1 \ldots i_k}|P'|$. \\

But $r'_{i_1 \ldots i_k}
   \le r^{\beta}_{i_1 \ldots i_k}$, therefore $d (x',y') \le r^{\beta}_{i_1 \ldots i_k}|P'| \le \left ( {\dfrac {d (x, y)} {\rho}}\right )^{\beta} |P'|$.\\

 b) There are   $ i_1 \ldots i_k i_{k+1}$ and $j_1 \ldots  j_l j_{l+1}$, such that  $x \in  P_{ i_1 \ldots i_k}\setminus P_{
 i_1 \ldots, i_k i_{k+1}}$,  $y \in  P_{j_1  \ldots,  j_l}\setminus P_{j_1 \ldots j_l j_{l+1}}$ and $P_{ i_1 \ldots i_k i_{k+1}} \bigcap  P_{j_1 \ldots
 j_l j_{l+1}} =S_{ i_1 \ldots i_k} (A)$,  where $A$ is some vertex of $P$.\\

In this case $d (x, y) \le \left
 \{ r_{i_1 \ldots i_k} +r_{j_1 \ldots j_l}\right\} |P|$. \\
  By Lemma  \ref{ga_xA},  $d (x, A) \ge  r_{i_1 \ldots i_k} \rho$ and $d (A, y) \ge {r_{j_1 \ldots j_l}} \rho$. \\ Therefore,  by Lemma \ref{ga_xy},  $d (x, y) \ge \rho\cdot { \sin {(\al/2)}} {( r_{i_1 \ldots i_k} + r_{j_1 \ldots j_l} )}$ , thus $$
  { (r_{i_1 \ldots i_k} + r_{j_1 \ldots j_l})} \rho\cdot {\sin{(\al/2)}} \le d (x, y) \le {(r_{i_1 \ldots i_k} +
 r_{j_1 \ldots j_l})} |P|.$$
 Similarly, for the system $\eS'$ we have
 $${( r'_{i_1 \ldots i_k} + r'_{j_1 \ldots j_l})} {\rho'}\cdot{\sin{({\al}'/2)}} \le d (x', y') \le {( r'_{i_1 \ldots i_k} + r'_{j_1 \ldots j_l})} |P'|.$$
 Suppose $ r_{i_1 \ldots i_k} \ge  r_{j_1
 \ldots j_l}$. Then, ${(r_{i_1 \ldots i_k})}{\rho}\cdot {\sin{(\al/2)}} \le d (x, y) \le 2 {( r_{i_1 \ldots i_k})} |P|$. \\

So, $d (x', y') \le 2  {(
 r_{i_1 \ldots i_k})}^{\beta} |P'| \le 2  \left ( \dfrac {d (x, y)} {\rho\cdot\sin {(\al/2)}}\right)^{\beta}|P'|.$ $\blacksquare$\\

{\vspace{.1in}
\noindent Mary Samuel\\
Department of Mathematics\\ 
Bharata Mata College, Kochi, India\\
marysamuel2000@gmail.com

\vspace{.1in}\noindent Andrei Tetenov\\
Gorno-Altaisk State University 
and\\
Novosibirsk State University, Russia\\
atet@mail.ru\\

\vspace{.1in}\noindent Dmitry Vaulin\\
Gorno-Altaisk State University\\
Gorno-Altaisk, Russia\\
d\_warrant@mail.ru }

\end{document}